\documentclass[10pt]{article}
\usepackage{amsfonts}
\usepackage{amsmath}
\usepackage{amscd,amssymb,amsthm}
\usepackage{graphics}
\usepackage{graphicx}

\newtheorem{theorem}{Theorem}[section]

\newtheorem{definition}[theorem]{Definition}
\newtheorem{example}[theorem]{Example}

\newtheorem{proposition}[theorem]{Proposition}
\newtheorem{remark}[theorem]{Remark}

\begin{document}

\title{Analysis of degenerate Chenciner bifurcation }
\author{G. Moza\thanks{
Department of Mathematics, Politehnica University of Timisoara, Romania;
email: gheorghe.moza@upt.ro}, S. Lugojan, L. Ciurdariu }
\date{}

\maketitle

\begin{abstract}
Degenerate Chenciner bifurcation in generic discrete-time dynamical systems
is studied in this work. While the non-degenerate Chenciner bifurcation can
be described by 2 bifurcation diagrams, the degeneracy we studied in this
work gives rise to 32 different bifurcation diagrams.
\end{abstract}


\section{Introduction}

The Chenciner bifurcation in a non-degenerate framework was initially
studied in \cite{che2} and \cite{che3}, and more recently in \cite{Ku98}, \cite%
{shi} and \cite{deng}. A description of the bifurcation can also be found in 
\cite{arrow}. The bifurcation arises in smooth discrete-time dynamical
systems (i.e. maps) as a generalization of the Neimark-Sacker bifurcation.
Discrete dynamical systems arise typically in 1) continuous-time dynamical
systems by discretizing the time \cite{bis} (particularly in numerical
simulations) or as return maps (such as Poincar\'{e} maps) defined by the
intersections of the systems' flows with some surfaces transversal to the
flows \cite{dum}, \cite{lib1}, \cite{mar}, \cite{men}, and 2) practical
applications when a phenomenon cannot be observed continuously in time but
at given moments of time \cite{cul}, \cite{jos}, \cite{mul}.

Bio-mathematical models based on discrete-time dynamical systems have long
been used in biology, particularly in population dynamics. For example,
methods to construct discrete dynamic models of gene regulatory networks
from experimental data sets are presented in \cite{kri}, while in \cite{fer}%
, a study on the dynamics of a mapping governing the time evolution of a
population spread over a one-dimensional lattice is reported. A
discrete-time epidemic model on heterogeneous networks has been recently
studied in \cite{wan}. The authors claim that the dynamics of endemic
equilibrium is more complex for discrete epidemic model compared to the
model based on differential equations. A discrete model to describe
epileptic seizures in a neural network of excitatory and inhibitory neurons
is proposed in \cite{vol}. In \cite{hui} a stroboscopic (discrete) map is
used to describe the existence of a infection-free periodic solution in an
epidemic model, while a prey--predator discrete model is studied in \cite%
{wan2}.

We aim to study in this work what happens to the Chenciner bifurcation when
it becomes \textit{degenerate} with respect to the transformation of
parameters, that is, when the transformation is not regular at $\left(
0,0\right).$ This form of degeneracy is also known as \textit{%
non-transversality}. A study concerning degenerate Hopf--Neimarck--Sacker
bifurcations can be found in \cite{bro1}.

The paper is organized in three main sections, following the introductory
section 1. In section 2 we analyze in a general setup the existence and
stability of equilibrium points and invariant closed curves when the
Chenciner bifurcation becomes degenerate. Section 3 describes the existence
of bifurcation curves and their dynamics in the parametric plane $\left(
\alpha _{1},\alpha _{2}\right) .$ We obtained a remarkable result in this
section, namely, two of the bifurcation curves coincide in their linear,
quadratic and cubic terms, and start to become two separate curves only with
terms of order four. Section 4 presents a wealth of bifurcation diagrams
which arise in the Chenciner bifurcation due to the degeneracy, which do not
exist in the non-degenerate framework.

\section{Analysis of degenerate Chenciner bifurcation}

To introduce the bifurcation, consider a difference equation of the form

\begin{equation}
x_{n+1}=f\left( x_{n},\alpha \right)  \label{ecn}
\end{equation}%
where $x_{n}\in 
%TCIMACRO{\U{211d} }%
%BeginExpansion
\mathbb{R}
%EndExpansion
^{2},$ $n\in 
%TCIMACRO{\U{2115} }%
%BeginExpansion
\mathbb{N}
%EndExpansion
,$ $\alpha =\left( \alpha _{1},\alpha _{2}\right) \in 
%TCIMACRO{\U{211d} }%
%BeginExpansion
\mathbb{R}
%EndExpansion
^{2}$ and $f$ is a smooth function of class $C^{r}$ with $r\geq 2.$ In order
to avoid indices, we often write \eqref{ecn} in the form 
\begin{equation}
x\longmapsto f\left( x,\alpha \right)  \label{ecs}
\end{equation}%
or $\tilde{x}=f\left( x,\alpha \right) .$ In complex coordinates, \eqref{ecs}
can be written as 
\begin{equation}
z\longmapsto \mu \left( \alpha \right) z+g\left( z,\bar{z},\alpha \right)
\label{ecz}
\end{equation}%
where $\mu $ and $g$ are smooth functions of their arguments, $\mu \left(
\alpha \right) =r\left( \alpha \right) e^{i\theta \left( \alpha \right) }$
where $r\left( 0\right) =1$ and $\theta \left( 0\right) =\theta _{0}.$ The
function $g$ can be written as a Taylor series 
\[
g\left( z,\bar{z},\alpha \right) =\sum_{k+l\geq 2}g_{kl}\left( \alpha
\right) z^{k}\bar{z}^{l} 
\]%
where $g_{kl}\left( \alpha \right)$ are smooth complex-valued functions.
Equation \eqref{ecz} can be further transformed to 
\begin{eqnarray}
w &\longmapsto &\left( r\left( \alpha \right) e^{i\theta \left( \alpha
\right) }+a_{1}\left( \alpha \right) w\bar{w}+a_{2}\left( \alpha \right)
w^{2}\bar{w}^{2}\right) w+O\left( \left\vert w\right\vert ^{6}\right)
\label{ecw} \\
&=&\left( r\left( \alpha \right) +b_{1}\left( \alpha \right) w\bar{w}%
+b_{2}\left( \alpha \right) w^{2}\bar{w}^{2}\right) we^{i\theta \left(
\alpha \right) }+O\left( \left\vert w\right\vert ^{6}\right)  \nonumber
\end{eqnarray}%
where $b_{k}\left( \alpha \right) =a_{k}\left( \alpha \right) e^{-i\theta
\left( \alpha \right) },$ $k=1,2.$ Denote by 
\begin{equation}
\beta _{1}\left( \alpha \right) =r\left( \alpha \right) -1\text{ and }\beta
_{2}\left( \alpha \right) =Re\left( b_{1}(\alpha )\right) .  \label{b12}
\end{equation}%
In polar coordinates \eqref{ecw} $\allowbreak $becomes

\begin{equation}
\left\{ 
\begin{array}{cc}
\rho _{n+1}= & \rho _{n}\left( 1+\beta _{1}\left( \alpha \right) +\beta
_{2}\left( \alpha \right) \rho _{n}^{2}+\allowbreak L_{2}\left( \alpha
\right) \rho _{n}^{4}\right) +\rho _{n}O\left( \rho _{n}^{6}\right) \\ 
\varphi _{n+1}= & \varphi _{n}+\theta \left( \alpha \right) +\rho
_{n}^{2}\left( \frac{Im\left( b_{1}\left( \alpha \right) \right) }{\beta
_{1}\left( \alpha \right) +1}+O\left( \rho _{n},\alpha \right) \allowbreak
\right)\ \ \ \ \ \ \ \ \ \text{ }%
\end{array}%
\right.  \label{rofi}
\end{equation}%
where $L_{2}\left( \alpha \right) =\frac{Im^{2}\left( b_{1}\left( \alpha
\right) \right) +2\left( 1+\beta _{1}\left( \alpha \right) \right) Re\left(
b_{2}\left( \alpha \right) \right) }{2\left( \beta _{1}\left( \alpha \right)
+1\right) }.$

A bifurcation in \eqref{rofi} which satisfies $r\left( 0\right) =1$ and $%
Re\left( b_{1}(0)\right) =0$ but $L_{2}\left( 0\right) \neq 0,$ is known as
the \textit{Chenciner bifurcation} (or generalized Neimark-Sacker
bifurcation). The Neimark-Sacker bifurcation arising in maps corresponds to
Hopf bifurcation \cite{Ku98}, \cite{sil}, \cite{waw} in continuous-time
dynamical systems, which is also related to fold-Hopf \cite{tig1}, \cite%
{tig2}, \cite{tig3} and zero-Hopf bifurcations \cite{lib2}.

\begin{remark}\label{re0} 
Discrete dynamical systems present complex behavior with chaos arising even in 1--dimensional (1D) maps \cite{fow}, while in 1D continuous-time dynamical systems chaos does not appear. In general, all objects existing in continuous-time dynamical systems, such as orbits, steady states, invariant sets, limit cycles, attractors, homoclinic, heteroclinic orbits and so on, have their correspondences in discrete-time dynamical systems.
\end{remark}

It follows from $\beta _{1}\left( 0\right) =0$ that 
\[
L_{2}\left( 0\right) =\frac{1}{2}\left( Im^{2}\left( b_{1}\left( 0\right)
\right) +2Re\left( b_{2}\left( 0\right) \right) \right) . 
\]

When the transformation of parameters 
\begin{equation}
\left( \alpha _{1},\alpha _{2}\right) \longmapsto \left( \beta _{1}\left(
\alpha \right) ,\beta _{2}\left( \alpha \right) \right)  \label{trp}
\end{equation}%
is regular at $\left( 0,0\right) ,$ then \eqref{rofi} can be put in a
simpler form. This is the \textit{non-degenerate} Chenciner bifurcation (see
e.g. \cite{Ku98} and \cite{shi}).

In this work we aim to approach the Chenciner bifurcation in the case when
the transformation (\ref{trp}) is not regular at $(0,0).$ The degeneracy
prevents us from using $\beta _{1,2}$ as new parameters. Instead, we have to
keep and work with the initial parameters $\alpha _{1,2}$ in the polar form %
\eqref{rofi}.

The truncated form of the $\rho -$map of \eqref{rofi} is 
\begin{equation}
\rho _{n+1}=\rho _{n}\left( 1+\beta _{1}\left( \alpha \right) +\beta
_{2}\left( \alpha \right) \rho _{n}^{2}+\allowbreak L_{2}\left( \alpha
\right) \rho _{n}^{4}\right) .  \label{tru}
\end{equation}%
The $\varphi -$map of \eqref{rofi} describes a rotation by an angle
depending on $\alpha $ and $\rho ,$ and can be approximated by 
\begin{equation}
\varphi _{n+1}=\varphi _{n}+\theta \left( \alpha \right) .  \label{truf}
\end{equation}%
Assume $0<\theta \left( 0\right) <\pi .$ The system to be studied in this
work is \eqref{tru}-\eqref{truf}, which is called the \textit{truncated
normal form} of the system \eqref{ecw}. The main equation shaping the
dynamics of the system \eqref{tru}-\eqref{truf} is the $\rho -$map %
\eqref{tru}. It is independent from the $\varphi -$map and can be studied
separately.

The $\rho -$map \eqref{tru} defines a one-dimensional dynamical system and
this will be further studied. It has the fixed point $\rho =0,$ for all
values of $\alpha ,$ which corresponds to the fixed point $O\left(
0,0\right) $ in the full system \eqref{tru}-\eqref{truf}. A positive nonzero
fixed point of the $\rho -$map \eqref{tru} corresponds to a closed invariant
curve (which is a circle) in the system \eqref{tru}-\eqref{truf}.

Note that $sign\left( L_{2}\left( \alpha \right) \right) =sign\left(
L_{0}\right) $ when $\left\vert \alpha \right\vert =\sqrt{\alpha
_{1}^{2}+\alpha _{2}^{2}}$ is sufficiently small because $L_{2}\left( \alpha
\right) =L_{0}\left( 1+O\left( \left\vert \alpha \right\vert \right) \right) 
$ and $L_{0}\neq 0.$ Throughout the work, $O\left( \left\vert \alpha
\right\vert ^{n}\right) $ for $n\geq 1$ denotes a series with real
coefficients $c_{ij}$ of the form 
\[
O\left( \left\vert \alpha \right\vert ^{n}\right) =\sum_{i+j\geq
n}c_{ij}\alpha _{1}^{i}\alpha _{2}^{j}. 
\]

\begin{proposition}
\label{pr1} The fixed point $O$ is (linearly) stable if $\beta _{1}\left(
\alpha \right) <0$ and unstable if $\beta _{1}\left( \alpha \right) >0,$ for
all values $\alpha $ with $\left\vert \alpha \right\vert $ sufficiently
small. On the bifurcation curve $\beta _{1}\left( \alpha \right) =0,$ $O$ is
(nonlinearly) stable if $\beta _{2}\left( \alpha \right) <0$ and unstable if 
$\beta _{2}\left( \alpha \right) >0,$ when $\left\vert \alpha \right\vert $
is sufficiently small. At $\alpha =0,$ $O$ is (nonlinearly) stable if $%
L_{0}<0$ and unstable if $L_{0}>0.$
\end{proposition}

PROOF. From $\left. \frac{d\rho _{n+1}}{d\rho _{n}}\right\vert _{\rho
=0}=1+\beta _{1}\left( \alpha \right) $ the proof follows when $\beta
_{1}\left( \alpha \right) \neq 0.$ When $\beta _{1}\left( \alpha \right) =0,$
the fixed point $O$ becomes nonhyperbolic and cannot be studied as in the
first case. Along the curve $\beta _{1}\left( \alpha \right) =0,$ the map %
\eqref{tru} reads 
\begin{equation}
\rho _{n+1}-\rho _{n}=\rho _{n}^{3}\left( \beta _{2}\left( \alpha \right)
+\allowbreak L_{2}\left( \alpha \right) \rho _{n}^{2}\right) .  \label{sm1}
\end{equation}%
By definition, $\rho _{n}\geq 0$ for all $n\geq 0.$ To study the stability
of $\rho =0,$ we consider orbits satisfying \eqref{sm1} and starting at
points $\rho _{0}$ sufficiently close to $\rho =0.$ For such orbits, when $%
\beta _{2}\left( \alpha \right) <0$ and $L_{0}<0$ we get $0<\rho _{n+1}<\rho
_{n},$ respectively, when $\beta _{2}\left( \alpha \right) >0$ and $L_{0}>0$
we have $0<\rho _{n}<\rho _{n+1},$ for all $\left\vert \alpha \right\vert $
sufficiently small. Thus, $\rho =0$ is stable in the first case and unstable
in the second.

Denote by $S_{n}=\beta _{2}\left( \alpha \right) +\allowbreak L_{2}\left(
\alpha \right) \rho _{n}^{2},$ for all $\left\vert \alpha \right\vert $
sufficiently small. When $\beta _{2}\left( \alpha \right) <0$ and $L_{0}>0,$
choose $\rho _{0}>0$ sufficiently small such as $S_{0}=\beta _{2}\left(
\alpha \right) +L_{2}\left( \alpha \right) \rho _{0}^{2}<0,$ that is, $%
\left\vert \rho _{0}\right\vert <\sqrt{\frac{-\beta _{2}\left( \alpha
\right) }{L_{2}\left( \alpha \right) }};$ $sign\left( L_{2}\left( \alpha
\right) \right) =sign\left( L_{0}\right) >0.$ It leads to $0<\rho _{1}<\rho
_{0}$ and $S_{1}<S_{0}$ because $\rho _{1}-\rho _{0}=\rho _{0}^{3}S_{0}<0$
and $S_{1}-S_{0}=\allowbreak L_{2}\left( \alpha \right) \left( \rho
_{1}^{2}-\rho _{0}^{2}\right) <0.$ Let us show by induction that $\rho
_{n+1}<\rho _{n}$ for all $n\geq 1.$ Indeed, assuming $\rho _{n+1}<\rho _{n}$
in \eqref{sm1} we obtain $S_{n}=\beta _{2}+\allowbreak L_{2}\rho _{n}^{2}<0$
and $S_{n+1}-S_{n}=L_{2}\left( \rho _{n+1}^{2}-\rho _{n}^{2}\right) <0$ for
all $\left\vert \alpha \right\vert $ is sufficiently small and $\rho _{n}>0.$
This yields $S_{n+1}<S_{n}<0$ which, in turn, implies $\rho _{n+2}-\rho
_{n+1}=\rho _{n+1}^{3}S_{n+1}<0$ and the conclusion follows.

Therefore, $0<\rho _{n+1}<\rho _{n}$ whenever $\beta _{2}\left( \alpha
\right) <0$ and $L_{0}\neq 0.$ One can similarly show that $0<\rho _{n}<\rho
_{n+1}$ when $\beta _{2}\left( \alpha \right) >0$ and $L_{0}\neq 0.$ These
results imply that $\rho =0$ is (nonlinearly) stable if $\beta _{2}\left(
\alpha \right) <0$ and unstable if $\beta _{2}\left( \alpha \right) >0$
whenever $L_{0}\neq 0.$

At $\alpha =0,$ \eqref{tru} yields $\rho _{n+1}-\rho _{n}=\allowbreak
L_{0}\rho _{n}^{5}\left( 1+O\left( \left\vert \alpha \right\vert \right)
\right) $ and the conclusion follows. $\blacksquare $

\bigskip

For each positive nonzero fixed point of the one-dimensional $\rho -$map %
\eqref{tru}, it corresponds an invariant curve (a circle) in the truncated
two-dimensional map \eqref{tru}-\eqref{truf}. The fixed points of \eqref{tru}
are solutions of the equation%
\begin{equation}
L_{2}\left( \alpha \right) y^{2}+\beta _{2}\left( \alpha \right) y+\beta
_{1}\left( \alpha \right) =0  \label{y1y2}
\end{equation}%
where $y=\rho _{n}^{2}.$ Denote by $\Delta \left( \alpha \right) =\beta
_{2}^{2}\left( \alpha \right) -4\beta _{1}\left( \alpha \right) L_{2}\left(
\alpha \right) ,$ respectively, $y_{1}=\frac{1}{2L_{2}}\left( \sqrt{\Delta }%
-\beta _{2}\right) $ and $y_{2}=-\frac{1}{2L_{2}}\left( \sqrt{\Delta }+\beta
_{2}\right) $ the two roots of \eqref{y1y2}, whenever they exist as real
numbers.

\begin{theorem}
\label{th1} The following assertions are true.

1) When $\Delta \left( \alpha \right) <0$ for all $\left\vert \alpha
\right\vert $ sufficiently small, the system \eqref{tru}-\eqref{truf} has no
invariant circles.

2) When $\Delta \left( \alpha \right) >0$ for all $\left\vert \alpha
\right\vert $ sufficiently small, the system \eqref{tru}-\eqref{truf} has:

\qquad a) one invariant unstable circle $\rho _{n}=\sqrt{y_{1}}$ if $L_{0}>0$ and $%
\beta _{1}\left( \alpha \right) <0;$

\qquad b) one invariant stable circle $\rho _{n}=\sqrt{y_{2}}$ if $L_{0}<0$ and $\beta
_{1}\left( \alpha \right) >0;$

\qquad c) two invariant circles, $\rho _{n}=\sqrt{y_{1}}$ unstable and $\rho _{n}=%
\sqrt{y_{2}}$ stable, if $L_{0}>0,$ $\beta _{1}\left( \alpha \right) >0,$ $%
\beta _{2}\left( \alpha \right) <0$ or $L_{0}<0,$ $\beta _{1}\left( \alpha
\right) <0,$ $\beta _{2}\left( \alpha \right) >0;$ in addition, $y_{1}<y_{2}$
if $L_{0}<0$ and $y_{2}<y_{1}$ if $L_{0}>0;$

\qquad d) no invariant circles if $L_{0}>0,$ $\beta _{1}\left( \alpha \right) >0,$ $%
\beta _{2}\left( \alpha \right) >0$ or $L_{0}<0,$ $\beta _{1}\left( \alpha
\right) <0,$ $\beta _{2}\left( \alpha \right) <0.$

3) On the bifurcation curve $\Delta \left( \alpha \right) =0,$ the system %
\eqref{tru}-\eqref{truf} has one invariant unstable circle $\rho _{n}=\sqrt{y_{1}}$
for all $L_{0}\neq 0.$ Moreover, if $L_{0}<0,$ the invariant circle is stable from the
exterior and unstable from the interior, while if $L_{0}>0$ is vice versa.

4) When $\beta _{1}\left( \alpha \right) =0,$ the system \eqref{tru}-%
\eqref{truf} has one invariant circle $\rho _{n}=\sqrt{-\frac{\beta _{2}\left( \alpha
\right) }{L_{0}}}$ whenever $L_{0}\beta _{2}\left( \alpha \right) <0.$ It is
stable if $L_{0}<0$ and $\beta _{2}\left( \alpha \right) >0,$ respectively,
unstable if $L_{0}>0$ and $\beta _{2}\left( \alpha \right) <0.$
\end{theorem}

PROOF. 1) - 2). Equation \eqref{y1y2} has no real roots when $\Delta \left(
\alpha \right) <0.$ Thus, the system \eqref{tru}-\eqref{truf} has no
invariant closed orbits (circles) in this case. Two different roots $y_{1,2}$
of \eqref{y1y2} exist iff $\Delta \left( \alpha \right) >0.$ They may lead
correspondingly to two circles in the truncated system \eqref{tru}-%
\eqref{truf}. The first circle $\rho _{n}=\sqrt{y_{1}}$ is unstable while
the second $\rho _{n}=\sqrt{y_{2}}$ is stable whenever they exist, i.e. when 
$y_{1,2}>0,$ because 
\begin{equation}
\left. \frac{d\rho _{n+1}}{d\rho _{n}}\right\vert _{\sqrt{y_{1}}}=1+2y_{1}%
\sqrt{\Delta }>1\text{ and }\left. \frac{d\rho _{n+1}}{d\rho _{n}}%
\right\vert _{\sqrt{y_{2}}}=1-2y_{2}\sqrt{\Delta }<1.  \label{der}
\end{equation}%
The items a)-d) follow from studying the signs of $y_{1}$ and $y_{2}$ from $%
y_{1}y_{2}=\frac{\beta _{1}\left( \alpha \right) }{L_{0}}$ and $y_{1}+y_{2}=-%
\frac{\beta _{2}\left( \alpha \right) }{L_{0}}.$

3) When $\Delta \left( \alpha \right) =0$ we can write \eqref{tru} in the
form 
\[
\rho _{n+1}-\rho _{n}=L_{0}\rho _{n}\left( \rho _{n}^{2}-\frac{\beta
_{2}\left( \alpha \right) }{2L_{2}\left( \alpha \right) }\right) ^{2}\left(
1+O\left( \left\vert \alpha \right\vert \right) \right) , 
\]%
which implies that $\rho _{n+1}\geq \rho _{n}$ when $L_{0}>0$ and $\rho
_{n+1}\leq \rho _{n}$ when $L_{0}<0,$ for all $\left\vert \alpha \right\vert 
$ sufficiently small. Thus, if $\Delta(\alpha) =0$ the circle $\rho _{n}=%
\sqrt{y_{1}},$ with $y_{1}=\frac{-\beta _{2}\left( \alpha \right) }{2L_{0}},$
is unstable for all $L_{0}\neq 0,$ whenever it exists, that is, on $%
L_{0}\beta _{2}\left( \alpha \right) <0.$ It is unstable because, when $%
L_{0}<0,$ an orbit starting at $\rho _{0}<\sqrt{y_{1}}$ will depart from the
circle since $\rho _{n+1}\leq \rho _{n}.$ The orbit through $\rho _{0}$ will
approach the origin $\rho =0.$ In fact, when $L_{0}<0,$ the orbits starting
from exterior points $\rho _{0}>\sqrt{y_{1}}$ to the circle $\rho _{n}=\sqrt{%
y_{1}}$ will tend to the circle, while orbits starting from interior points $%
\rho _{0}<\sqrt{y_{1}}$ will depart from the circle. Thus, the circle is
stable from the exterior and unstable from the interior. When $L_{0}>0,$ the
scenario is opposite: the orbits starting at $\rho _{0}>\sqrt{y_{1}}$ will
depart from the circle due to $\rho _{n+1}\geq \rho _{n},$ while orbits
starting at $\rho _{0}<\sqrt{y_{1}}$ will approach the circle. Thus, the
circle is unstable from the exterior and stable from the interior. Figure %
\ref{f3} presents the generic phase portraits of this case.

4) When $\beta _{1}\left( \alpha \right) =0$ a circle may exist when $%
L_{0}\beta _{2}\left( \alpha \right) <0.$ Indeed, when $L_{0}<0$ and $\beta
_{2}\left( \alpha \right) >0,$ one obtain $y_{1}=0$ and $y_{2}=-\frac{\beta
_{2}}{L_{0}}>0,$ while for $L_{0}>0$ and $\beta _{2}\left( \alpha \right) <0$
we have $y_{1}=-\frac{\beta _{2}}{L_{0}}>0$ and $y_{2}=0.$ Thus, in the both
cases, the circle is $\rho _{n}=\sqrt{-\frac{\beta _{2}}{L_{0}}}.$ It
follows from \eqref{der} that the circle is stable if $L_{0}<0$ and $\beta
_{2}\left( \alpha \right) >0$ and unstable if $L_{0}>0$ and $\beta
_{2}\left( \alpha \right) <0.$ In addition, from Proposition \ref{pr1}, $O$
is unstable if $\beta _{2}\left( \alpha \right) >0$ and stable if $\beta
_{2}\left( \alpha \right) <0.$ $\blacksquare $

\section{Bifurcation curves}

Write in the following the smooth functions $\beta _{1,2}\left( \alpha
\right) $ as Taylor series at $\left( 0,0\right) $ in the form

\begin{equation}
\beta _{1}\left( \alpha \right) =\sum_{i+j=1}^{p}a_{ij}\alpha _{1}^{i}\alpha
_{2}^{j}+O\left( \left\vert \alpha \right\vert ^{p+1}\right) \text{ and }%
\beta _{2}\left( \alpha \right) =\sum_{i+j=1}^{q}b_{ij}\alpha _{1}^{i}\alpha
_{2}^{j}+O\left( \left\vert \alpha \right\vert ^{q+1}\right)  \label{b1b2g}
\end{equation}%
for some $p,q\geq 1.$

The transformation \eqref{trp} is not regular at $\left( 0,0\right) $ and,
thus, the Chenciner bifurcation is \textit{degenerate}, if and only if $%
\left. \frac{\partial \beta _{1}}{\partial \alpha _{1}}\frac{\partial \beta
_{2}}{\partial \alpha _{2}}-\frac{\partial \beta _{1}}{\partial \alpha _{2}}%
\frac{\partial \beta _{2}}{\partial \alpha _{1}}\right\vert _{\alpha =0}=0,$
that is, 
\begin{equation}
a_{10}b_{01}-a_{01}b_{10}=0.  \label{ns}
\end{equation}%
The behavior of the system \eqref{ecw} is unknown in this case because the
bifurcation diagram obtained for the non-degenerate Chenciner bifurcation 
\cite{Ku98} is not valid anymore when \eqref{ns} is satisfied. In fact, we
will show in this work that the new constraint \eqref{ns} changes
significantly the bifurcation diagrams of the system \eqref{ecw}, without
changing its generic phase portraits. While in the non-degenerate framework,
one or two bifurcation diagrams were sufficient to describe the behavior, in
the case of the degeneracy a large number of bifurcation diagrams are needed
for this purpose.

Denote by $B_{1,2}$ and $C$ the following sets of points in $%
%TCIMACRO{\U{211d} }%
%BeginExpansion
\mathbb{R}
%EndExpansion
^{2}$ 
\begin{equation}
B_{1,2}=\left\{ \left( \alpha _{1},\alpha _{2}\right) \in 
%TCIMACRO{\U{211d} }%
%BeginExpansion
\mathbb{R}
%EndExpansion
^{2},\beta _{1,2}\left( \alpha \right) =0,\left\vert \alpha \right\vert
<\varepsilon \right\}  \label{b12g}
\end{equation}%
and 
\begin{equation}
C=\left\{ \left( \alpha _{1},\alpha _{2}\right) \in 
%TCIMACRO{\U{211d} }%
%BeginExpansion
\mathbb{R}
%EndExpansion
^{2},\Delta \left( \alpha \right) =0,\left\vert \alpha \right\vert
<\varepsilon \right\}  \label{cg}
\end{equation}%
for some $\varepsilon >0$ sufficiently small. The following result is
crucial for determining bifurcation diagrams, thus, for understanding the
behavior of the truncated system \eqref{tru}-\eqref{truf}. It describes the
existence and relative positions one to another in the parametric plane $%
\left( \alpha _{1},\alpha _{2}\right) $ of bifurcation curves.

\begin{definition}
We say that two curves of the form $\alpha _{2}=\sum_{i=1}^{p}c_{i}\alpha
_{1}^{i}+c_{p+1}\alpha _{1}^{p+1}+O\left( \alpha _{1}^{p+2}\right) $ and $%
\alpha _{2}=\sum_{i=1}^{p}b_{i}\alpha _{1}^{i}+b_{p+1}\alpha
_{1}^{p+1}+O\left( \alpha _{1}^{p+2}\right) $ coincide in their first $p$
terms if and only if $c_{i}=b_{i}$ for all $1\leq i\leq p$ and $c_{p+1}\neq
b_{p+1}.$
\end{definition}

\begin{theorem}
\label{th2} Assume the degeneracy condition \eqref{ns} holds true and $%
a_{01}b_{01}\neq 0.$ Then, the sets $B_{1,2}$ and $C$ are smooth curves of
the form $\alpha _{2}=c_{1}\alpha _{1}+O\left( \alpha _{1}^{2}\right) ,$ $%
c_{1}=-\frac{a_{10}}{a_{01}},$ tangent to the same line $a_{10}\alpha
_{1}+a_{01}\alpha _{2}=0.$ Moreover, $C$ and $B_{1}$ coincide up to their
cubic terms in $\alpha _{1}$ and may become two separated curves in the
parametric plane only with terms in $\alpha _{1}^{4}$ or higher.
\end{theorem}

PROOF. When $a_{01}\neq 0,$ the curve $\beta _{1}\left( \alpha \right) =0,$
denoted by $B_{1},$ is well-defined and unique by Implicit Function Theorem
(IFT). It can be written for $\left\vert \alpha _{1}\right\vert $ small
enough in the form 
\begin{equation}
B_{1}=\left\{ \left( \alpha _{1},\alpha _{2}\right)\in\mathbb{R}%
^{2},\alpha_{2}=c_{1}\alpha _{1}+c_{2}\alpha _{1}^{2}+c_{3}\alpha
_{1}^{3}+c_{4}\alpha _{1}^{4}+O\left( \alpha _{1}^{5}\right) \right\}
\label{b10}
\end{equation}%
where $c_{1}=-\frac{a_{10}}{a_{01}}.$ Similarly, the curve $\beta _{2}\left(
\alpha \right) =0$ exists and is unique provided that $b_{01}\neq 0,$ and
can be expressed in the form

\begin{equation}
B_{2}=\left\{ \left( \alpha _{1},\alpha _{2}\right) \in 
%TCIMACRO{\U{211d} }%
%BeginExpansion
\mathbb{R}
%EndExpansion
^{2},\alpha _{2}=d_{1}\alpha _{1}+d_{2}\alpha _{1}^{2}+O\left( \alpha
_{1}^{3}\right) \right\}  \label{b20}
\end{equation}%
where $d_{1}=-\frac{b_{10}}{b_{01}}=c_{1}.$ From $d_{1}=c_{1},$ the curve $%
B_{2}$ is tangent at $O$ to the same line $a_{10}\alpha _{1}+a_{01}\alpha
_{2}=0$ as $B_{1}.$ An important quantity in the further analysis is $\Delta
\left( \alpha \right) =\beta _{2}^{2}\left( \alpha \right) -4\beta
_{1}\left( \alpha \right) L_{2}\left( \alpha \right) ,$ which in this case
reads 
\begin{equation}
\Delta \left( \alpha \right) =-4L_{0}\left( a_{10}\alpha _{1}+a_{01}\alpha
_{2}\right) +O\left( \left\vert \alpha \right\vert ^{2}\right)
\label{delta0}
\end{equation}%
where $L_{0}\neq 0.$ Using the same argument based on the Implicit Function
Theorem, the bifurcation curve $\Delta \left( \alpha \right) =0$ is
well-defined and unique, and can be written for $\left\vert \alpha
_{1}\right\vert $ small enough in the form%
\begin{equation}
C=\left\{ \left( \alpha _{1},\alpha _{2}\right) \in 
%TCIMACRO{\U{211d} }%
%BeginExpansion
\mathbb{R}
%EndExpansion
^{2},\alpha _{2}=c_{1}\alpha _{1}+m_{2}\alpha _{1}^{2}+m_{3}\alpha
_{1}^{3}+m_{4}\alpha _{1}^{4}+O\left( \alpha _{1}^{5}\right) \right\} .
\label{c10}
\end{equation}%
We notice that $C$ is also tangent at $O$ to the same line as $B_{1}$ and $%
B_{2},$ having the same coefficient $c_{1}$ of its lowest term. The three
curves $C$ and $B_{1,2}$ coincide when they are approximated by their linear
terms, being given by $\alpha _{2}=c_{1}\alpha _{1}.$ This property is due
to the degeneracy condition \eqref{ns}. In order to draw bifurcation
diagrams we need to determine $C$ and $B_{1,2}$ as distinct curves in the
parametric plane $\alpha_1O\alpha_2.$ Thus, we need to find the second order
terms containing $\alpha _{1}^{2}$ in the expressions of the curves.

To this end, write $\beta _{1}\left( \alpha \right) =a_{10}\alpha
_{1}+a_{01}\alpha _{2}+a_{20}\alpha _{1}^{2}+a_{11}\alpha _{1}\alpha
_{2}+a_{02}\alpha _{2}^{2}+O\left( \left\vert \alpha \right\vert ^{3}\right)
,$ that is, $p=2$ in \eqref{b1b2g}. We can determine $c_{2}$ by deriving
with respect to $\alpha _{1}$ at $0$ the expression $\beta _{1}\left( \alpha
\right) =0,$ considering that $\alpha _{2}$ is a function of $\alpha _{1}$
such that $\alpha _{2}=0$ if $\alpha _{1}=0$ (from IFT). Alternatively, we
can substitute for $\alpha _{2}=c_{1}\alpha _{1}+c_{2}\alpha _{1}^{2}$ in $%
\beta _{1}\left( \alpha \right) =0$ and find $c_{2}$ by nullifying the
coefficient of $\alpha _{1}^{2}$ in 
\begin{equation}
\beta _{1}\left( \alpha \right) =\left( a_{10}+c_{1}a_{01}\right) \alpha
_{1}+\left( a_{02}c_{1}^{2}+a_{11}c_{1}+a_{20}+c_{2}a_{01}\right)
\allowbreak \alpha _{1}^{2}+O\left( \alpha _{1}^{3}\right) .  \label{b1s}
\end{equation}%
Similarly, the curve $B_{2}$ has the form $\alpha _{2}=c_{1}\alpha
_{1}+d_{2}\alpha _{1}^{2}+O\left( \alpha _{1}^{3}\right) .$ After
calculations we obtain 
\begin{equation}
c_{2}=-\frac{1}{a_{01}}\left( a_{20}+c_{1}^{2}a_{02}+c_{1}a_{11}\right) 
\text{ and }d_{2}=-\frac{1}{b_{01}}\left(
b_{20}+c_{1}^{2}b_{02}+c_{1}b_{11}\right) .  \label{c2d2}
\end{equation}%
It is clear that $d_{2}\neq c_{2}$ in general. Assume $d_{2}\neq 0.$ Thus,
the two curves $B_{1}$ and $B_{2}$ are different starting with terms of
order $\alpha _{1}^{2}.$ Having known $c_{2}$ and $d_{2},$ the curves $%
B_{1,2}$ can be determined as two separated curves in the parametric plane.

Let us determine the quadratic term $m_{2}\alpha _{1}^{2}$ in the
bifurcation curve $C$ given by \eqref{c10}. To this end, we need $p=2,$ $q=1$
in \eqref{b1b2g} and $L_{2}\left( \alpha \right) =L_{0}+l_{10}\alpha
_{1}+l_{01}\alpha _{2}+O\left( \left\vert \alpha \right\vert ^{2}\right) .$
After calculations and using $c_{1}=-\frac{a_{10}}{a_{01}}=-\frac{b_{10}}{%
b_{01}}$ we obtain $m_{2}=c_{2}.$ This implies that the bifurcation curves $%
C $ and $B_{1}$ coincide up to their quadratic terms. Thus, we have to
determine their cubic terms.

Using $p=3$ in \eqref{b1b2g} we find 
\begin{equation}
c_{3}=-\frac{1}{a_{01}}\left(
a_{30}+c_{1}^{2}a_{12}+c_{1}^{3}a_{03}+c_{1}a_{21}+c_{2}a_{11}+2c_{1}c_{2}a_{02}\right) .
\label{c3}
\end{equation}%
In order to determine cubic terms in $\Delta$ and find $m_{3}$ for the curve 
$C,$ we need to use $p=3$ and $q=2$ in \eqref{b1b2g}, and $L_{2}\left(
\alpha \right) =L_{0}+\sum_{i+j=1}^{2}a_{ij}\alpha _{1}^{i}\alpha
_{2}^{j}+O\left( \left\vert \alpha \right\vert ^{3}\right) .$ Surprisingly,
we find $m_{3}=c_{3}.$ Therefore, the curve $\Delta \left( \alpha \right) =0$
coincides to the curve $\beta _{1}\left( \alpha \right) =0$ up their cubic
terms, namely, they both have the form $\alpha _{2}=c_{1}\alpha
_{1}+c_{2}\alpha _{1}^{2}+c_{3}\alpha _{1}^{3}+O\left( \alpha
_{1}^{4}\right) .$ At this point we may think that a global property is
hidden in the relationship of the curves $C$ and $B_{1},$ which would make
them to coincide fully when $\left\vert \alpha \right\vert $ is sufficiently
small. However, this is a misleading conclusion because, looking at the term 
$\beta _{2}^{2}\left( \alpha \right) $ in the expression of $\Delta \left(
\alpha \right) =\beta _{2}^{2}\left( \alpha \right) -4\beta _{1}\left(
\alpha \right) L_{2}\left( \alpha \right) ,$ we notice that $\beta
_{2}\left( \alpha \right) =\left( b_{10}+c_{1}b_{01}\right) \alpha
_{1}+O\left( \alpha _{1}^{2}\right) $ when $\left( \alpha _{1},\alpha
_{2}\right) \in C.$ But from the degeneracy condition \eqref{ns} we have $%
b_{10}+c_{1}b_{01}=0,$ which means that $\beta _{2}^{2}\left( \alpha \right)
=O\left( \alpha _{1}^{4}\right) ,$ that is, the term $\beta _{2}^{2}\left(
\alpha \right) $ contributes to the expression of $\Delta \left( \alpha
\right) $ only with terms of order $4$ or higher. This is the reason why the
two curves $C$ and $B_{1}$ coincide up to cubic terms. The exact calculation
of the next terms in $C$ and $B_1$ is given in the followings. Using $p=4$
in \eqref{b1b2g} we find

\[
c_{4}=-\frac{1}{a_{01}}\left( a_{04}c_{1}^{4}+a_{13}c_{1}^{3}+\left(
a_{22}+3c_{2}a_{03}\right) c_{1}^{2}+\sigma _{1}\allowbreak c_{1}+\sigma
_{2}\right) , 
\]
where $\sigma _{1}=a_{31}+2c_{2}a_{12}+2c_{3}a_{02}\allowbreak $ and $\sigma
_{2}=a_{02}c_{2}^{2}+a_{21}c_{2}+a_{40}+c_{3}a_{11}.$

Further, in order to determine $m_{4}$ from the expression of $C$ \eqref{c10}%
, we use $p=4$ and $q=2$ in \eqref{b1b2g}, respectively $L_{2}\left( \alpha
\right) =L_{0}+\sum_{i+j=1}^{3}a_{ij}\alpha _{1}^{i}\alpha _{2}^{j}+O\left(
\left\vert \alpha \right\vert ^{4}\right) .$ With these ingredients, we
obtain

\begin{equation}
m_{4}=c_{4}+\allowbreak \frac{\gamma ^{2}}{4L_{0}a_{01}}  \label{m4}
\end{equation}%
where $\gamma =b_{02}c_{1}^{2}+b_{11}c_{1}+b_{20}+c_{2}b_{01}.$ Assume $%
m_{4}\neq 0$ and $c_{4}\neq 0.$ The equality \eqref{m4} is important because
it determines the relativ positions of the bifurcation curves $C$ and $B_{1}$
in the parametric plane, whenever $\gamma \neq 0$ and $L_{0}a_{01}c_{4}\neq
0.$ $\blacksquare $

\section{Bifurcation diagrams}

Assume the linear parts of $\beta _{1,2}\left( \alpha \right) $ have nonzero
coefficients, that is, $a_{10}b_{01}a_{01}b_{10}\neq 0$ such that $%
a_{10}b_{01}=a_{01}b_{10}.$ Thus, the three bifurcation curves are
well-defined and unique when $\left\vert \alpha \right\vert $ is
sufficiently small. Assume also $d_{2}m_{4}c_{4}\gamma L_{0}\neq 0.$ Two
main cases arise for obtaining bifurcation diagrams, corresponding to $%
L_{0}a_{01}>0$ and $L_{0}a_{01}<0.$

The bifurcation curve $C$ splits the parametric plane into two regions $%
R_{1,2}.$ Denote by $R_{1}$ the region where $\alpha _{2}>0$ and $R_{2}$
with $\alpha _{2}<0.$

\begin{figure}[h!]
\begin{center}
\includegraphics[width=0.45\textwidth]{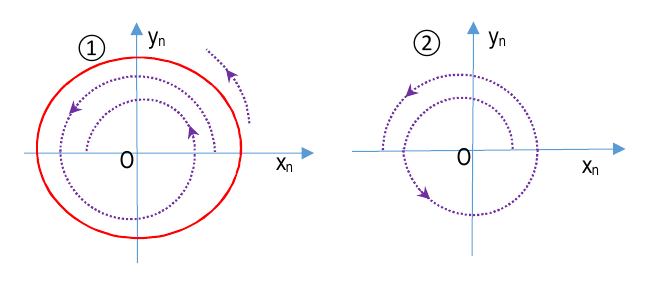} \includegraphics[width=0.45%
\textwidth]{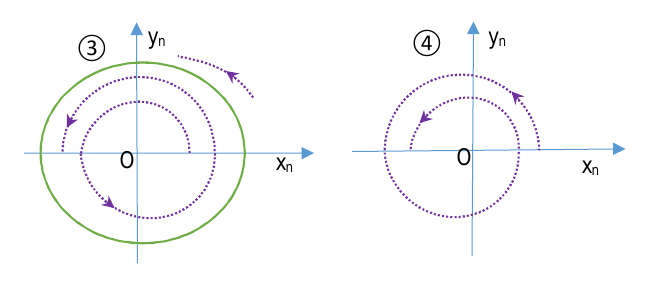} \includegraphics[width=0.45\textwidth]{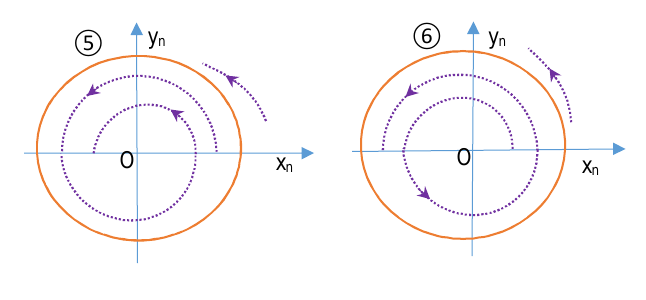} %
\includegraphics[width=0.45\textwidth]{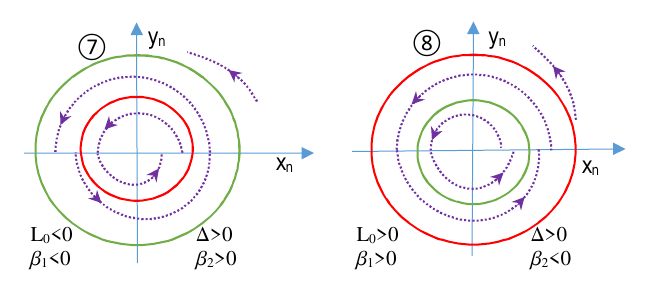}
\end{center}
\caption{Generic phase portraits for $\protect\theta_0>0.$}
\label{f1}
\end{figure}

\begin{remark}
\label{re1} Figure \ref{f1} presents generic phase portraits corresponding
to different regions of the bifurcation diagrams, including the phase
portraits on the bifurcation curves defined by $\Delta \left( \alpha \right)
=0,$ respectively, $\beta _{1}\left( \alpha \right) =0.$ We summarize in
Table 1 the correspondence between $\Delta ,$ $\beta _{1,2},$ $L_{0}$ and
the generic phase portraits, respectively, different regions from
bifurcation diagrams. When $\beta _{1,2}\left( \alpha \right) =0$ then $%
\alpha =0.$
\end{remark}

\begin{center}
Tabel 1. \textit{The regions in the parametric plane defined by $%
\Delta(\alpha),$ $\beta_{1,2}(\alpha)$ and $L_{0}.$}

\bigskip

\begin{tabular}{lllll|llllll}
$\Delta \left( \alpha \right) $ & $L_{0}$ & $\beta _{1}\left( \alpha \right) 
$ & $\beta _{2}\left( \alpha \right) $ & $Region$ &  & $\Delta \left( \alpha
\right) $ & $L_{0}$ & $\beta _{1}\left( \alpha \right) $ & $\beta _{2}\left(
\alpha \right) $ & $Region$ \\ \hline
$+$ & $+$ & $+$ & $+$ & $2$ &  & $0$ & $-$ & $-$ & $-$ & $4$ \\ 
$+$ & $-$ & $-$ & $-$ & $4$ &  & $0$ & $-$ & $-$ & $+$ & $5$ \\ 
$+$ & $+$ & $-$ & $\pm ,0$ & $1$ &  & $0$ & $+$ & $+$ & $-$ & $6$ \\ 
$+$ & $-$ & $+$ & $\pm ,0$ & $3$ &  & $0$ & $+$ & $0$ & $0$ & $2$ \\ 
$+$ & $-$ & $-$ & $+$ & $7$ &  & $0$ & $-$ & $0$ & $0$ & $4$ \\ 
$+$ & $+$ & $+$ & $-$ & $8$ &  & $+$ & $-$ & $0$ & $+$ & $3$ \\ 
$-$ & $+$ & $+$ & $\pm ,0$ & $2$ &  & $+$ & $-$ & $0$ & $-$ & $4$ \\ 
$-$ & $-$ & $-$ & $\pm ,0$ & $4$ &  & $+$ & $+$ & $0$ & $-$ & $1$ \\ 
$0$ & $+$ & $+$ & $+$ & $2$ &  & $+$ & $+$ & $0$ & $+$ & $2$%
\end{tabular}
\end{center}

\begin{remark}\label{re1c}
The \textit{circles} presented in Figure \ref{f1} are closed invariant curves which resemble \textit{limit cycles} from continuous-time dynamical systems. However, there exist several fundamental differences between the two types of curves. Firstly, the circles from Figure \ref{f1} (the closed curves) are filled discretely with points while limit cycles are filled continuously. A limit cycle is made up of a single periodic orbit, while on a closed invariant curve it might exist an infinite number of discrete orbits, periodic or not. 
\end{remark}

\begin{figure}[h!]
\begin{center}
\includegraphics[width=70mm, height=55mm]{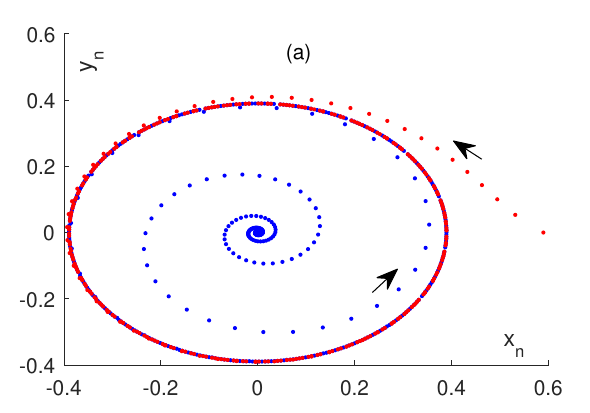} %
\includegraphics[width=70mm, height=55mm]{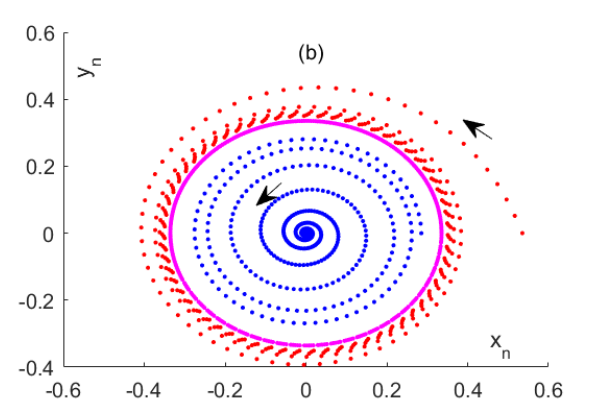}
\end{center}
\caption{Numerical simulations for the map (\protect\ref{ex1}) with $%
\protect\beta_{1,2}(\protect\alpha)=\protect\alpha_1+\protect\alpha_2$ and $\protect\alpha_1=\protect\alpha_2=0.01$ in case (a), respectively, $\beta_{1,2}$ given by (\ref{ex1a}), $a_1=0.5$ and $a_2=-0.512$ in case (b).}
\label{f12}
\end{figure}

\begin{example} We illustrate numerically in Figure \ref{f12} the existence of closed
invariant curves in the particular two-dimensional map given in polar
coordinates by 
\begin{equation}
\rho _{n+1}=\rho _{n}\left( 1+\beta _{1}\left( \alpha \right) +\beta
_{2}\left( \alpha \right) \rho _{n}^{2}-\rho _{n}^{4}\right) \text{ and }%
\varphi _{n+1}=\varphi _{n}+\theta _{0},  \label{ex1}
\end{equation}%
where $\left\vert \alpha \right\vert $ is sufficiently small and $\theta
_{0}=0.1;$ $L_{0}=-1.$ In Figure \ref{f12} (a) we consider $\beta
_{1,2}\left( \alpha \right) =\alpha _{1}+\alpha _{2},$ which yield $\Delta >0
$ on $\alpha _{1}+\alpha _{2}>0,$ while in Figure \ref{f12} (b)  
\begin{equation}
\beta_{1}=\alpha _{1}+\alpha _{2}-\left( \alpha _{1}+\alpha _{2}\right)
^{2}\  \text{and}\ \beta _{2}=2\left( \alpha _{1}+\alpha
_{2}\right) +\alpha _{1}^{2}.\label{ex1a}
\end{equation}%
Thus, in the second case, the curve $(C)$
defined by $\Delta \left( \alpha \right) =0$ reads $4\alpha _{2}\left(
\alpha _{1}^{2}+1\right) =\allowbreak -\alpha _{1}\left( 4\alpha
_{1}^{2}+\alpha _{1}^{3}+4\right) ,$ which in the lowest terms becomes $%
\alpha _{2}=\allowbreak -\alpha _{1}-\frac{1}{4}\alpha _{1}^{4}+O\left(
\alpha _{1}^{5}\right) .$ One can check that $\left. \beta _{1}\right\vert
_{C}<0$ and $\left. \beta _{2}\right\vert _{C}>0.$

We notice that in the both cases the transformation (\ref{trp}) is not
regular at $\left( 0,0\right),$ thus, the Chenciner bifurcation is
degenerate in (\ref{ex1}). The discrete orbits in Figure \ref{f12} (a) tend
to a circle (the invariant stable closed curve), while in Figure \ref{f12}
(b), some orbits approach the circle (those starting at exterior points to
the circle), while others depart from it (the interior orbits). 

Figure \ref{f12} (a) illustrates the generic phase portrait 3 while (b) the
portrait 5 from Figure \ref{f1}. Notice the circle in Figure \ref{f12} (b)
is stable from the exterior and unstable from the interior. In Figure \ref%
{f12} we used the Cartesian coordinates $(x_{n},y_{n}),$ with $x_{n}=\rho
_{n}\cos \varphi _{n}$ and $y_{n}=\rho _{n}\sin \varphi _{n}.$
\end{example}

\begin{figure}[h!]
\begin{center}
\includegraphics[width=0.24\textwidth]{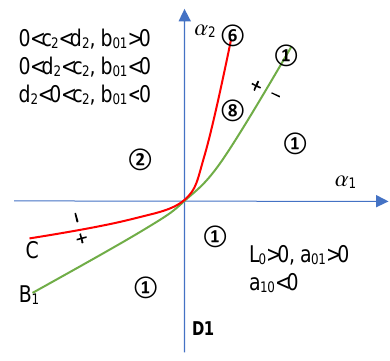} \includegraphics[width=0.24%
\textwidth]{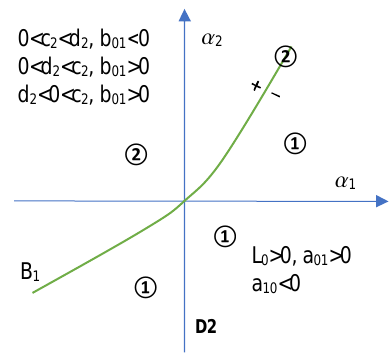} \includegraphics[width=0.24\textwidth]{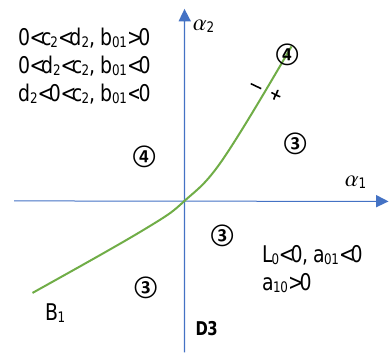} %
\includegraphics[width=0.24\textwidth]{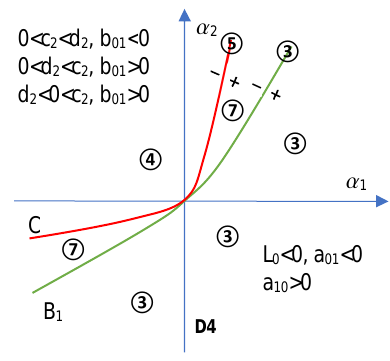} \includegraphics[width=0.24%
\textwidth]{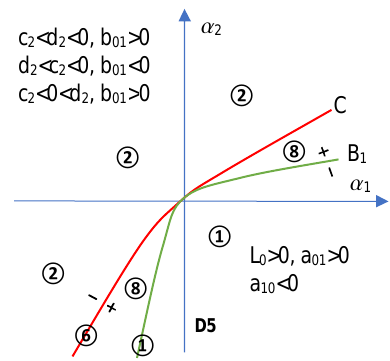} \includegraphics[width=0.24\textwidth]{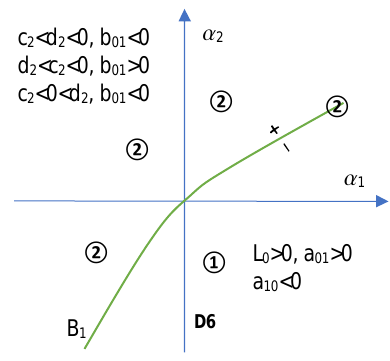} %
\includegraphics[width=0.24\textwidth]{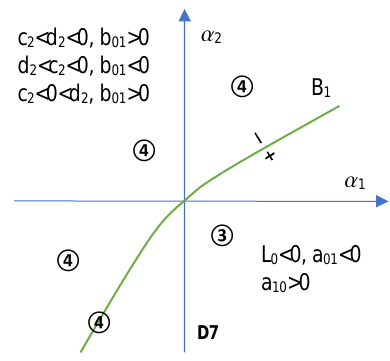} \includegraphics[width=0.24%
\textwidth]{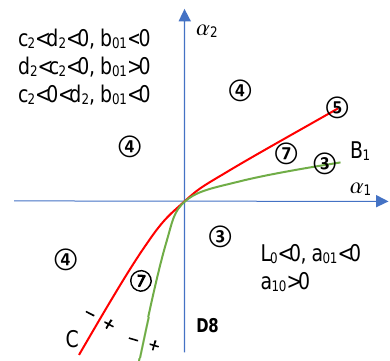}
\end{center}
\caption{Bifurcation diagrams corresponding to $L_{0}a_{01}>0$ and $c_1>0.$ }
\label{f2}
\end{figure}

\begin{figure}[h!]
\begin{center}
\includegraphics[width=0.24\textwidth]{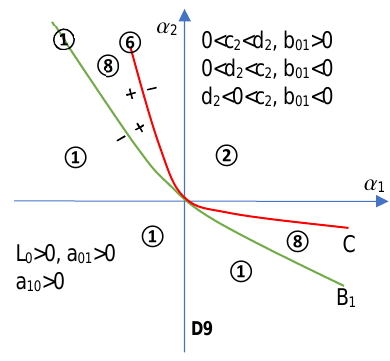} \includegraphics[width=0.24%
\textwidth]{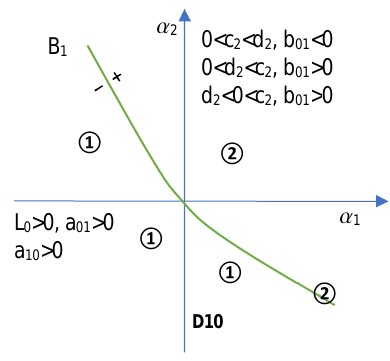} \includegraphics[width=0.24\textwidth]{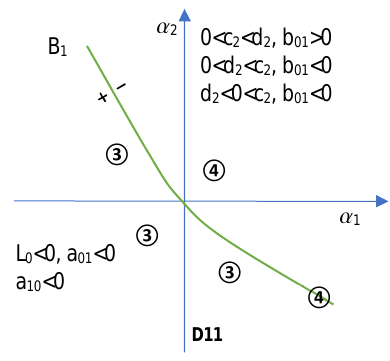} %
\includegraphics[width=0.24\textwidth]{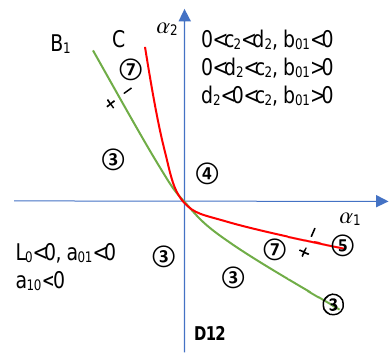} \includegraphics[width=0.24%
\textwidth]{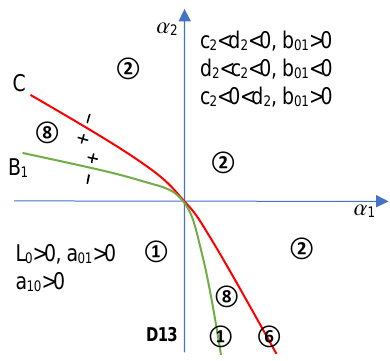} \includegraphics[width=0.24\textwidth]{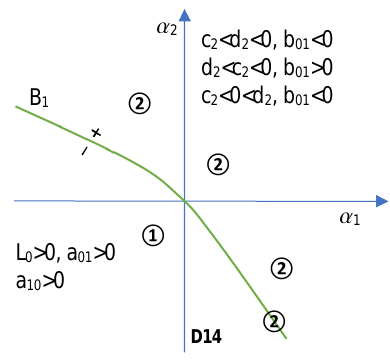} %
\includegraphics[width=0.24\textwidth]{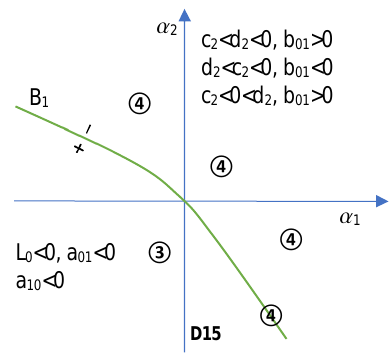} \includegraphics[width=0.24%
\textwidth]{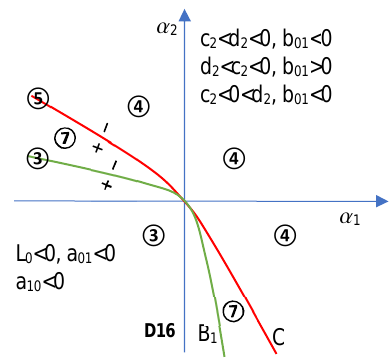}
\end{center}
\caption{Bifurcation diagrams corresponding to $L_{0}a_{01}>0$ and $c_1<0.$ }
\label{f3}
\end{figure}

\textbf{Case 1.} When $L_{0}a_{01}>0,$ then $m_{4}>c_{4}$ by \eqref{m4}.
From \eqref{delta0} it follows that $\Delta \left( \alpha \right) <0$ on $%
R_{1}$ and $\Delta \left( \alpha \right) >0$ on $R_{2},$ when $\left\vert
\alpha \right\vert $ is sufficiently small. Figures \ref{f2}-\ref{f3}
contain the bifurcation diagrams of this case. Notice that, the bifurcation
diagrams can be obtained for different conditions. For example, $D1$
corresponds to $0<c_2<d_2,$ $b_{01}>0$ or $0<d_2<c_2,$ $b_{01}<0$ or $%
d_2<0<c_2,$ $b_{01}<0.$

\begin{figure}[h!]
\begin{center}
\includegraphics[width=0.24\textwidth]{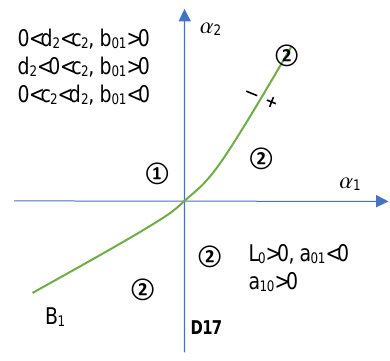} \includegraphics[width=0.24%
\textwidth]{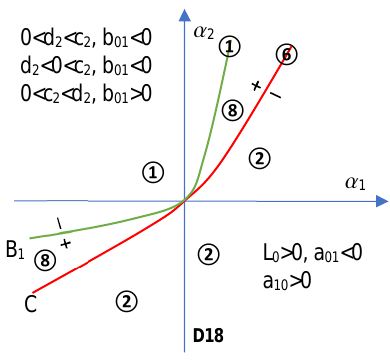} \includegraphics[width=0.24\textwidth]{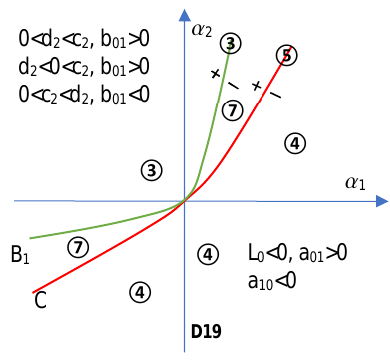} %
\includegraphics[width=0.24\textwidth]{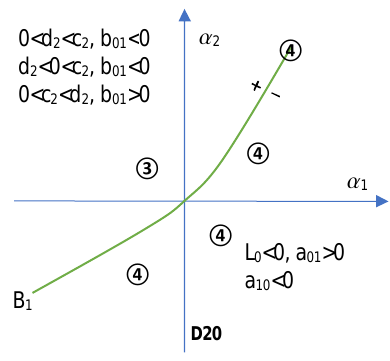} \includegraphics[width=0.24%
\textwidth]{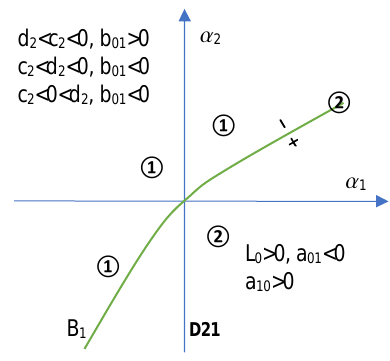} \includegraphics[width=0.24\textwidth]{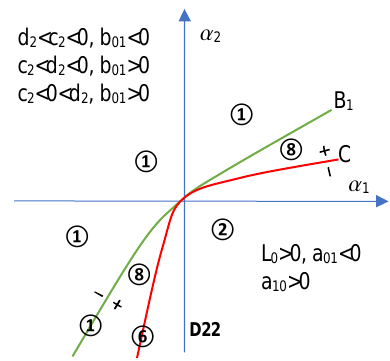} %
\includegraphics[width=0.24\textwidth]{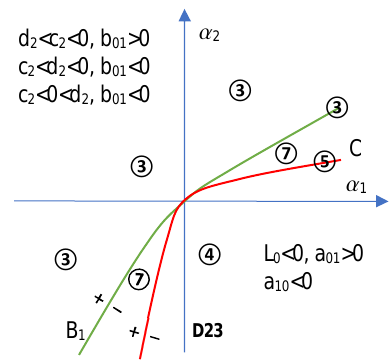} \includegraphics[width=0.24%
\textwidth]{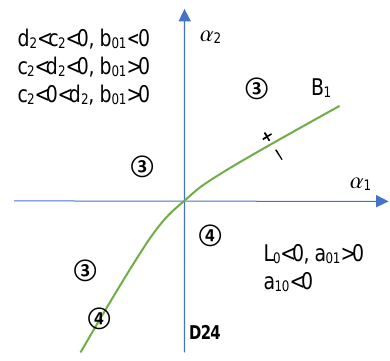}
\end{center}
\caption{Bifurcation diagrams corresponding to $L_{0}a_{01}<0$ and $c_1>0.$ }
\label{f4}
\end{figure}

\begin{figure}[h!]
\begin{center}
\includegraphics[width=0.24\textwidth]{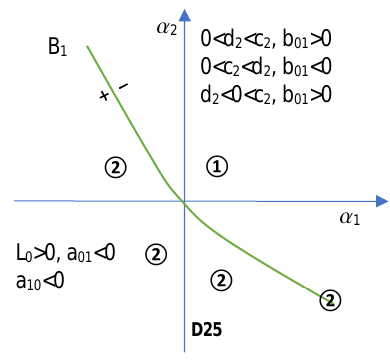} \includegraphics[width=0.24%
\textwidth]{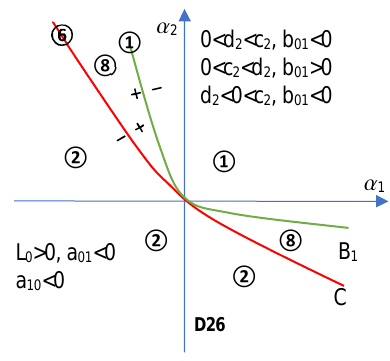} \includegraphics[width=0.24\textwidth]{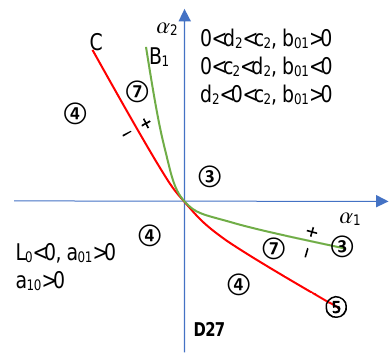} %
\includegraphics[width=0.24\textwidth]{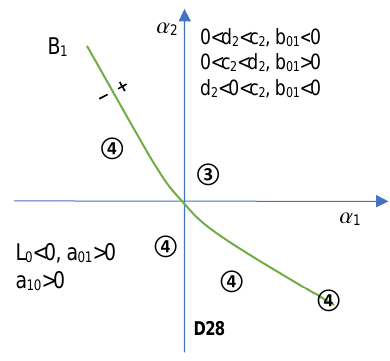} \includegraphics[width=0.24%
\textwidth]{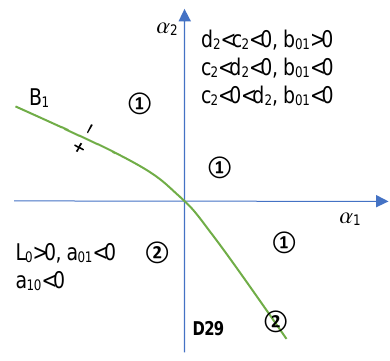} \includegraphics[width=0.24\textwidth]{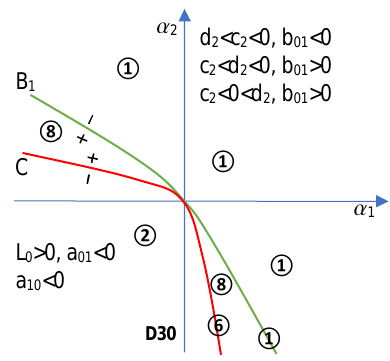} %
\includegraphics[width=0.24\textwidth]{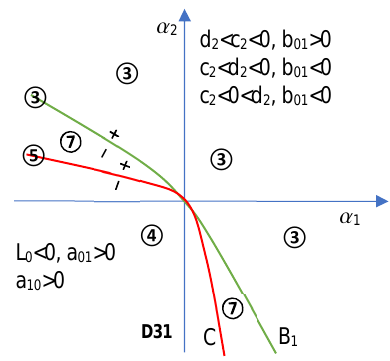} \includegraphics[width=0.24%
\textwidth]{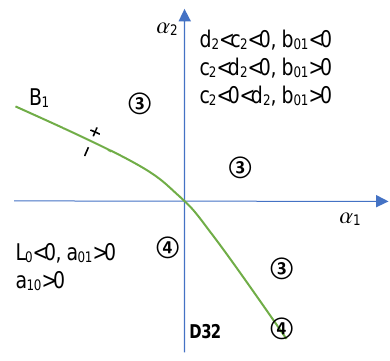}
\end{center}
\caption{Bifurcation diagrams corresponding to $L_{0}a_{01}<0$ and $c_1<0.$ }
\label{f5}
\end{figure}

\textbf{Case 2.} When $L_{0}a_{01}<0,$ then $m_{4}<c_{4}$ by \eqref{m4}. In
this case $\Delta \left( \alpha \right) >0$ on $R_{1}$ and $\Delta \left(
\alpha \right) <0$ on $R_{2},$ when $\left\vert \alpha \right\vert $ is
sufficiently small. We present in Figures \ref{f4}-\ref{f5} all bifurcation
diagrams corresponding to $L_{0}a_{01}<0$ and $c_1>0$ or $c_1<0.$

\begin{remark}
When $a_{10}=b_{10}=0$ and $a_{01}b_{01}\neq 0,$ one obtain similar
bifurcation diagrams to Cases 1-2. The tangent line $a_{10}\alpha
_{1}+a_{01}\alpha _{2}=0$ coincides to the horizontal $\alpha _{1}-$axis and
the bifurcation curves $C$ and $B_{1,2}$ become tangent to the $\alpha _{1}-$%
axis in this case.
\end{remark}

\section{Conclusions}

In this paper we studied the truncated normal form of Chenciner bifurcation
in a case of degeneracy, known also as the case of non-transversality. More
exactly, we assumed that the transformation of parameters used in
determining the normal form is not regular at $\alpha =0.$ While in the
non-degenerate framework the behavior of the system can be described by two
bifurcation diagrams, the degeneracy considered in this work gives rise to a
large number of bifurcation diagrams needed to describe the Chenciner
bifurcation.

We presented in this work bifurcation diagrams obtained in a general case,
namely, when the functions $\beta _{1}\left( \alpha \right) $ and $\beta
_{2}\left( \alpha \right) $ have both non-zero linear terms which satisfy
the degeneracy condition \eqref{ns}. We could also draw conclusions on $%
a_{10}=b_{10}=0.$ But what happens in other cases, such as $a_{10}=a_{01}=0$
or $b_{10}=b_{01}=0$ or $a_{10}=a_{01}=b_{10}=b_{01}=0,$ remains an open
problem. The study of the degenerate Chenciner bifurcation becomes more
difficult in these cases due to the involved bifurcation curves. For
example, when $\beta _{1}\left( \alpha \right) =\sum_{i+j=2}^{p}a_{ij}\alpha
_{1}^{i}\alpha _{2}^{j}+O\left( \left\vert \alpha \right\vert ^{p+1}\right)
, $ we need to find adequate procedures which establish when $\beta
_{1}\left( \alpha \right) =0$ represents a curve in the parametric plane.
This task is difficult because the Implicit Function Theorem cannot be
applied in general.

The results obtained in this work for the truncated normal form provides an
approximate description of the complicated bifurcation structure near a
generic Chenciner bifurcation \cite{che4}, \cite{bae}. More studies have to
be performed for the full normal form \eqref{rofi}.

\section{Acknowledgments}

This research was supported by Horizon2020-2017-RISE-777911 project.


\bigskip

\begin{thebibliography}{99}


%{Chenciner \textit{et al.}(1987)}

\bibitem{arrow} Arrowsmith, D. \& Place, C.
[1990], \textit{An Introduction to Dynamical Systems}, (Cambridge University
Press, Cambridge, 1990).

\bibitem{bae} Baesens C. \& MacKay, R.S. [2007]
``Resonances for weak coupling of the unfolding of a saddle-node periodic
orbit with an oscillator'', \textit{Nonlinearity} \textbf{20(5)}, 1283--1298.

\bibitem{bis} Biswas, M. \& Bairagi, N. [2020] ``On
the dynamic consistency of a two-species competitive discrete system with
toxicity'', \textit{Journal of Computational and Applied Mathematics} \textbf{%
363}, 145--155.

\bibitem{bro1} Broer, H. W., Holtman, S. J.,
Vegter, G. \& Vitolo, R. [2009] ``Geometry and dynamics of mildly degenerate
Hopf--Neimarck--Sacker families near resonance'', \textit{Nonlinearity} 
\textbf{22}, 2161--2200.

\bibitem{che2} Chenciner, A. [1985] ``Bifurcations de
points fixes elliptiques. II. Orbites periodiques et ensembles de Cantor
invariants'', \textit{Invent. Math.} \textbf{80}, 81--106.

\bibitem{che3} Chenciner, A. [1988] ``Bifurcations de points
fixes elliptiques. III. Orbites periodiques de petites periodes,'' \textit{%
Inst. Hautes Etudes Sci. Publ. Math.} \textbf{66}, 5-91.

\bibitem{che4} Chenciner, A., Gasull, A. \&
Llibre, J. [1987] ``Une description complete du portrait de phase d'un modele
d'elimination resonante'', \textit{C. R. Acad. Sci. Paris Ser. I Math.} 
\textbf{305 (13)}, 623--626.

\bibitem{cul} Cull, P., Flahive, M., Robson, R.
[2005] \textit{Difference Equations}, (Springer--Verlag New York, 2005).

\bibitem{deng} Deng, S. [2019] ``Bifurcations of a bouncing ball
dynamical system'', \textit{International Journal of Bifurcation and Chaos}, 
\textbf{29(14)}, 1950191.

\bibitem{dum} Dumortier, F., Llibre, J.,
Artes, J. C. [2006] \textit{Qualitative Theory of Planar Differential Systems%
}, (Springer--Verlag, Berlin, 2006).

\bibitem{fer} Fernandez, B. \& Tereshko, V.
[2000] ``Population dynamics in heterogeneous environments: A discrete
model'', \textit{International Journal of Bifurcation and Chaos}, \textbf{%
10(8)}, 1993--2000.

\bibitem{fow} Fowler, A. \& McGuinness, M.
[2020] \textit{One-dimensional Maps}, (Chaos, Springer, Cham).

\bibitem{hui} Hui, J. \& Zhu D. [2007] ``Dynamics
of SEIS epidemic models with varying population size'', \textit{International
Journal of Bifurcation and Chaos}, \textbf{17(05)}, 1513--1529.

\bibitem{jos} Jost, J. [2014] \textit{Mathematical Methods in
Biology and Neurobiology}, (Springer--Verlag London, 2014).

\bibitem{kri} Kriete, A. \& Roland Eils, R. [2014] 
\textit{Computational Systems Biology}, 2nd Edition, (Academic Press, 2014).

\bibitem{Ku98} Kuznetsov, Y.A. [2004] \textit{Elements of
Applied Bifurcation Theory}, 3rd Ed. (Springer--Verlag, 2004).

\bibitem{lib1} Llibre, J. and Teruel, A. E. [2014] 
\textit{Introduction to the Qualitative Theory of Differential Systems:
Planar, Symmetric and Continuous Piecewise Linear Systems},
(Springer--Verlag, 2014).

\bibitem{lib2} Llibre, J., Oliveira, R.D.S.,
\& Valls, C. [2015] ``On the integrability and the zero--Hopf bifurcation of
a Chen--Wang differential system'', \textit{Nonlinear Dynamics} \textbf{80},
353--361.

\bibitem{mar} Martelli, M. [2011] \textit{Introduction to
Discrete Dynamical Systems and Chaos}, (John Wiley \& Sons, 2011).

\bibitem{men} Mencinger, M., Fercec, B.,
Oliveira, R., Pagon, D. [2017] ``Cyclicity of some analytic maps'', \textit{%
Applied Mathematics and Computation} \textbf{295}, 114--125.

\bibitem{mul} Muller, J. \& Kuttler, C. [2015] 
\textit{Methods and Models in Mathematical Biology}, (Springer--Verlag
Berlin Heidelberg, 2015).

\bibitem{shi} Shilnikov, L.P., Shilnikov,
A.L., Turaev, D.V., Chua, L.O. [2001] \textit{Methods of Qualitative Theory
in Nonlinear Dynamics}, part 2, (World Scientific, 2001).

\bibitem{sil} Silva, V.B., Vieira, J.P.,
Leonel, E.D. [2018] ``A new application of the normal form description to a
N-dimensional dynamical systems attending the conditions of a Hopf
bifurcation'', \textit{Journal of Vibration Testing and System Dynamics} 
\textbf{2(3)}, 249--256.

\bibitem{tig1} Tigan, G. [2017] ``Degenerate with respect to
parameters fold-Hopf bifurcations'', \textit{Discrete and Continuous
Dynamical Systems--Series A} \textbf{37(4)}, 2115--2140.

\bibitem{tig2} Tigan, G., Llibre, J. \&
Ciurdariu, L. [2017] ``Degenerate fold-Hopf bifurcations in a Rossler-type
system'', \textit{International Journal of Bifurcation and Chaos} \textbf{%
27(5)}, 1--8.

\bibitem{tig3} Tigan, G. [2018] ``Analysis of degenerate
fold--Hopf bifurcation in a three--dimensional differential system'', \textit{%
Qualitative Theory of Dynamical Systems} \textbf{17(2)}, 387--402.

\bibitem{vol} Volk, D. [1999] ``Epileptic seizures in a discrete
model of neural networks of the brain'', \textit{International Journal of
Modern Physics}, \textbf{10(05)}, 815--821.

\bibitem{wan} Wang, X., Lu, J., Wang, Z., Li, Y.
[2020] ``Dynamics of discrete epidemic models on heterogeneous networks'', 
\textit{Physica A: Statistical Mechanics and its Applications}, \textbf{539}%
, 122991.

\bibitem{wan2} Wang, J. \& Feckan, M. [2020] ``Dynamics
of a discrete nonlinear prey--predator model'', \textit{International Journal
of Bifurcation and Chaos} \textbf{30(04)}, 2050055.

\bibitem{waw} Wawrzaszek, A. \& Krasinska, A.
[2019] ``Hopf bifurcations, periodic windows and intermittency in the
generalized Lorenz model'', \textit{International Journal of Bifurcation and
Chaos} \textbf{29(14)}, 1930042.
\end{thebibliography}
\end{document}